\documentclass[12pt]{amsart}
\usepackage[dvips]{graphicx}
\usepackage{epsfig}
\usepackage{amssymb}
\usepackage[usenames, dvipsnames]{color}
\usepackage{hyperref}
\usepackage{verbatim}
\usepackage[normalem]{ulem}

\theoremstyle{plain}
\newtheorem{theorem}{Theorem}[section]

\theoremstyle{definition}

\theoremstyle{remark}
\newtheorem{remark}[theorem]{Remark}

\def\lam{\lambda}

\def\th{\theta}

\def\e{\epsilon}

\def\p{\partial}

\def\f{\frac}

\def\al{\alpha}

\def\O{\Omega}

\def\q{\quad}
\def\qq{\qquad}

\def\be{\begin{equation}}
\def\ee{\end{equation}}
\def\bes{\begin{equation*}}
\def\ees{\end{equation*}}
\def\al{\begin{aligned}}
\def\eal{\end{aligned}}

\def\lab{\label}

\def\2O{\underline{\O}}

\def\lam{\lambda}
\def\sig{\sigma}

\numberwithin{equation}{section}

\makeatletter
\def\dashint{\operatorname%
{\,\,\text{\bf--}\kern-.98em\DOTSI\intop\ilimits@\!\!}}
\makeatother

\begin{document}


\title[type II model]{On a potential type II singular model of the Navier Stokes equations
}

\thanks{}

\author[Q. S. Zhang]{Qi S. Zhang}

\address[]{Department of mathematics, University of California, Riverside, CA 92521,
USA}

\email{qizhang@math.ucr.edu}



\begin{abstract} It is shown that a particular class of partial type I  or type II functions are not singularity models for the three dimensional Navier-Stokes equations with or  without smooth, compactly supported forces. This extends  part of the results for the axially symmetric case  to the full three dimension case.  Some connections to the recent literature are also discussed.

\end{abstract}
\maketitle


\section{Statement of result}

 In this paper,
  we study a potential singular model  of  the
  3 dimensional  Navier-Stokes equations:
\be
\label{nse}
\Delta v -  (v\cdot \nabla) v - \nabla P -\partial_t v =f, \quad \text{div} \, v=0, \quad \mathbf{R}^3 \times (0, \infty), \, v(x, 0)=v_0(x),
\ee where  $v=v(x, t)$ is the velocity and $v_0$ is the initial value; $P=P(x, t) $ is the pressure; $f=f(x, t)$ is the forcing term (or force).
 In the 1930s,  Leray \cite{Le2}, proved that the Cauchy
  problem of \eqref{nse} with $f=0$ has a solution in the energy space $(v, \nabla v) \in   (L^\infty_t L^2_x,
  L^2_{x t})$. However, as of this writing, in general it is not known if such solutions with zero force stay bounded or
  regular for all $t>0$ given regular initial values.

 However, most recently, on the eighth of September of 2026, OpenAi posted a paper \cite{Op} that draws wide attention. It states that there exists a smooth, compactly supported force $f$ such that \eqref{nse} with zero initial value $v_0=0$ has a smooth solution in the energy space for some time, but which becomes singular in finite time.
 While the validity of the result is being verified by people in the pertinent fields, it is clear that the paper already provides fresh ideas and valuable  information for the study of Navier Stokes equation. See the useful papers \cite{LR1}, \cite{LR2}  for further details and expositions on \cite{Op}.
 Curiously, the core part of the singular solution in \cite{Op} falls into the so called {\it partial type I} solutions defined earlier in the paper \cite{Zq} which was posted in April 2026.
 Recall that a solution to \eqref{nse} is of partial type I if, in the cylindrical system \eqref{ert3}, only the radial velocity $v_r$ satisfies the one-sided, critical bound
 \be
 \lab{vr>-c}
 v_r(x, t) \ge -\f{C}{\sqrt{T-t}}, \quad x \in \mathbf{R}^3, \q 0 \le t <T.
 \ee Here $C, T$ are positive constants. The concept is partly motivated by the type I regularity results for the axially symmetric Navier Stokes equations in \cite{CSTY1}, \cite{CSTY2} and \cite{KNSS}, where all components of the velocity is assumed to be bounded by $\f{C}{\sqrt{T-t}}$.

  In \cite{Zq}, it was proven that in the axially symmetric case (defined below \eqref{eqasns}), finite energy solutions are always regular provided the force is $0$ and the initial angular momentum density $r v_\theta |_{t=0}$ is bounded. Here $v_\theta$ is the angular speed in the cylindrical system. In fact the same regularity result is also true with smooth, compactly supported force. The proof is essentially the same. See also a related,  localized result with different proof in \cite{CIV1} and \cite{CIV2}.
Hence, if its result is true, there must be a very delicate blow up mechanism in the OpenAi paper, using bounded, smooth angular force.

In Theorem \ref{th1.1} below, it is shown that a particular partial type I blow up solution \eqref{sscene}, in the full 3 dimensional case, does not exist  with or without smooth, compactly supported force $f$ .
Based on this, in Theorem \ref{th1.2} below, under some extra conditions, it is shown that a type of genuinely  type II  ( or super critical) blow up solutions \eqref{sscene2} does not exist. This type of nonexistent blow up solutions  closely resembles that in \cite{Op}.  Since we have not finished reading the paper \cite{Op}, we do not pretend to have found a  counter example. See also Remark \ref{rmk2} below. The purpose is to foster discussion and eventual understanding. In any case, ruling in or ruling out potential singularity models is always an useful undertaking. For instance, the various forms of regularity criteria in the literature such as \cite{CKN:1} and the $L^{p, q}_{x, t, loc}$ ($3/p+ 2/q \le 1$) conditions by Ladyzhenskaya-Prodi-Serrin,  are, in essence, to rule out certain singularity models. In recent years, super critical criterions have appeared in \cite{Pa} (log supercritical axially symmetric case), \cite{Ta} (log supercritical 3 dimensional case), \cite{Ib} (super critical derivative condition) and  \cite{Se} (Euler type with extra conditions).

 In order to proceed, we also need  to use the cylindrical system. Let
 \[
 v = v_{r}e_{r} + v_{\th}e_{\th} + v_{3}e_{3}
 \]be the velocity in the cylindrical system with the standard basis $\{e_{r}, e_{\th}, e_{3}\}$, where for any $ x=(x_1,x_2,x_3)\in \mathbf{R}^3 $, $ r=\sqrt{x_1^2+x_2^2} \equiv |x'|$, $\theta= \tan^{-1} (x_2/x_1)$, and
\be\label{ert3}
e_r=(x_1/r, x_2/r, 0),\quad e_\th=(-x_2/r, x_1/r, 0),\quad e_3=(0,0,1).\ee
Then the Navier-Stokes equations become: for $r>0, x_3 \in \mathbf{R}, \th \in [0, 2 \pi]$,
\be
\lab{eqasns}
\begin{aligned}
\begin{cases}
   \big (\Delta-\frac{1}{r^2} \big )
v_r-(v_r \p_r + \f{1}{r} v_\th\p_\theta + v_3 \p_{x_3})v_r+\frac{(v_{\theta})^2}{r} -
 \f{2}{r^2} \p_\theta v_\th -\partial_r
P-\p_t  v_r=f_r,  \\
   \big   (\Delta-\frac{1}{r^2}  \big
)v_{\theta}-(v_r \p_r + \f{1}{r} v_\th \p_\theta + v_3 \p_{x_3} )v_{\theta}-\frac{v_{\theta} v_r}{r} + \f{2}{r^2} \p_\theta  v_r -  \f{1}{r} \p_\theta P-
\partial_t v_{\theta}=f_\th,\\
 \Delta v_3-(v_r \p_r + \f{1}{r} v_\th \p_\theta + v_3 \p_{x_3})v_3-\p_{x_3} P-\p_t v_3=f_3,\\
 \frac{1}{r} \p_r (rv_r)+ \f{1}{r} \p_\theta v_\th +\p_{x_3}
v_3=0,\\
v(x, 0)=v_0(x).
\end{cases}
\end{aligned}
\ee  In the above $f=f_r e_r + f_\th e_\th + f_3 e_3$.
If the unknown functions $P$, $f$, $v_{r}$, $v_{\th}$ and $v_{3}$ are independent of the azimuthal angle $\th$, we call the above axially symmetric Navier-Stokes equations, abbreviated as ASNS. In this case, all terms involving $\p_\th$ are $0$. The regularity problem for ASNS is still open at the time of this writing.  Although it was realized some years ago that the vortex stretching term is critical after a suitable change of variables \cite{LZ17}, building on \cite{CFZ}.

Before stating the main result, we need to introduce a few more notations.
For simplicity of presentation, we pick $[-1, 0)$ as the time interval in which the solutions  are smooth with $t=0$ being the potentially singular moment. The variable $\tau$ denotes $-t$. Following \cite{Op}, we will use the similarity variables $q, X, \eta$ defined as follows:
\be
\lab{xetaq}
X=\f{r^2}{2 q}, \q \eta=\f{x_3}{q^{0.5-h}}, \q \tau= q (1-\eta^2)=-t,
\ee where $h \in (0, 0.01)$ is a fixed number. See the beginning of Section 2 for more transparent relations between $X, \eta$, $q$ and the physical variables $x, t$. It can  be easily shown that $q=q(x_3, t)$ and
\be
\lab{taudenj0}
\tau + x^{2/(1-2h)}_3 \le q \le C (\tau + x^{2/(1-2h)}_3)
\ee and $(X, \eta, q) \in [0, \infty) \times (-1, 1) \times [0, \infty)$.

The following two theorems are the main results of the paper.
The potential singularity models consist of  principal terms which are axially symmetric plus  lower order terms depending on the angle $\th$. Their profiles are functions of the similarity variables $X, \eta$, the angle $\theta$ and time $t$. Note that $q$ is not an independent variable due to the 2nd and 3rd identity in \eqref{xetaq}, but it plays important roles through out.

\begin{theorem}
\lab{th1.1}
Let $v=v(x, t)=v_{r}e_{r} + v_{\th}e_{\th} + v_{3}e_{3}$ be a smooth Leray-Hopf solution of the  Navier-Stokes equations in $\mathbf{R^3} \times [-1, 0)$ with a force $f$ that is bounded, smooth and compactly supported in the spatial variables.
Suppose: for $t \in [-1, 0)$,

(a). the functions $E$, $\sqrt{X} \p_X E$, $\p_\eta E$ of variables $(X, \eta)$; $R$  and $U$ of $(X, \eta, \th, t)$  are smooth, bounded  and converge uniformly on compact sets of $\{ X>0 \}$ as $t \uparrow 0$;

(b). the functions
$E_2$, $\sqrt{X} \p_X E_2$, $\p_\eta E_2$  of $(X, \eta, \th, t)$
are smooth, bounded and converge uniformly to $0$ on compact sets of $\{ X>0 \}$ as $t \uparrow 0$; moreover, there are positive numbers $\delta$ and $C$ such that
$|E_2| \le C q^\delta$;

(c).
\[
\sup_{X>0, \eta \in [-1, 1], \th \in [0, 2\pi], t \in [0, 1)} |X|^{0.5^+} (|E(X, \eta)+ E_2(X, \eta, \th, t)|)<\infty
\]where $0.5^+$ is  a number strictly greater than $0.5$.

Then there does not exist a singularity scenario with the following asymptotic profiles
\be
\lab{sscene}
\al
&v_r = \frac{1}{q^{0.5}}  R(X, \eta, \th, t),\\
&v_\th = \frac{1}{q^{0.5+h}} \left(E(X, \eta) + E_2(X, \eta, \th, t)\right),\\
&v_3 = \frac{1}{q^{0.5+h}}  U(X, \eta, \th, t), \qq h \in (0, 0.01),
\eal
\ee such that
\be
\lab{thrvth}
\limsup_{t \to 0^-} |t|^h \Vert r v_\th(\cdot, t) \Vert_\infty>0.
\ee
\end{theorem}

Some elucidation of the theorem is given in the following remarks.

{ \remark
\lab{rmk0.1}
 In the above theorem, the condition on the radial velocity $v_r$ is critical under the natural scaling of the Navier Stokes equations. Therefore the solution is partial type I.
 Furthermore, in the next theorem, under mild extra conditions, the condition for $v_r$ becomes a super critical one.  The conditions on the angular and axial velocity $v_\th$ and $v_3$ are slightly super critical in both theorems. }

{ \remark
\lab{rmk0.2}
(1). The appearance of the function $\sqrt{X} \p_X E$ is due to the convention that $X=r^2/(2 q)$ and $E$ is  a smooth function of $\sqrt{X}$ or $r$ in the physical variable when $q \neq 0$. Note $q=0$ is equivalent to both $t=0$ and $x_3=0$.

(2). Condition (b) for $E_2$ is modeled on the perturbation term for the profile of the angular speed $v_\th$ in \cite{Op}. See 3.3 on p11 there for a description. The number $\delta$ here is chosen as $h/2$ there. The condition $|E_2| \le C q^\delta$ corresponds to the construction that $E_2$ contains a factor $q^\delta$.
 See Remark \ref{rmk2} below for more details.
 }

\begin{theorem}
\lab{th1.2}
Let $v=v(x, t)=v_{r}e_{r} + v_{\th}e_{\th} + v_{3}e_{3}$ be a smooth Leray-Hopf solution of the  Navier-Stokes equations in $\mathbf{R^3} \times [-1, 0)$ with a force $f$ that is bounded, smooth and compactly supported in the spatial variables.
Suppose: for $t \in [-1, 0)$,

(a). the functions $E$, $\sqrt{X} \p_X E$, $\p_\eta E$ of variables $(X, \eta)$; $R$  and $U$ of $(X, \eta, \th, t)$  are smooth, bounded  and converge uniformly on compact sets of $\{ X>0 \}$ as $t \uparrow 0$;

(b). the functions
$E_2$, $\sqrt{X} \p_X E_2$, $\p_\eta E_2$  of $(X, \eta, \th, t)$
are smooth, bounded and converge uniformly to $0$ on compact sets of $\{ X>0 \}$ as $t \uparrow 0$; moreover, there are positive numbers $\delta$ and $C$ such that
$|E_2| \le C q^\delta$;

(c).
\[
\sup_{X>0, \eta \in [-1, 1], \th \in [0, 2\pi], t \in [0, 1)} |X|^{0.5^+} (|E(X, \eta)+ E_2(X, \eta, \th, t)|)<\infty
\]where $0.5^+$ is  a number strictly greater than $0.5$;

(d). the functions $\f{1}{r} \p_\th E_2$, $\sqrt{X} \p_X U$, $\p_\eta U$ are bounded;

(e). $R=R(\cdot, 0)$ is not identically $0$ for $(X, \eta, \th) \in (0, \infty) \times (-1, 1) \times [0, 2 \pi]$.

Then there does not exist a singularity scenario with the following asymptotic profiles
\be
\lab{sscene2}
\al
&v_r = \frac{1}{q^{\sig}}  R(X, \eta, \th, t), \q \sig \ge 0.5,\\
&v_\th = \frac{1}{q^{0.5+h}} \left(E(X, \eta) + E_2(X, \eta, \th, t)\right),\\
&v_3 = \frac{1}{q^{0.5+h}}  U(X, \eta, \th, t), \qq h \in (0, 0.01),
\eal
\ee such that
\be
\lab{thrvth02}
\limsup_{t \to 0^-} |t|^h \Vert r v_\th(\cdot, t) \Vert_\infty>0.
\ee
\end{theorem}

Observe that Conditions (d) and (e) are the only extra ones comparing with those of Theorem \ref{th1.1}. The conditions in Theorem \ref{th1.2} are genuinely super critical or type II.

{\remark
\lab{rmk1}
 In the interesting papers \cite{CIV2} and \cite{IKX},  regularity of solutions to \eqref{nse} under the same anisotropic conditions as \cite{Op} were proven with additional conditions  such as analyticity of the force or rotationally invariant curl of the force respectively. The proof is based on the Euler scaling method in \cite{Se}, which was further developed in \cite{CIV1}. These authors actually replaced the asymptotically self-similar profiles by more general conditions on the three components of the velocity and their derivatives.

In contrast, we assume the asymptotically self-similar profiles \eqref{sscene2} hold, but make no additional assumption on the force. The compactly supported condition on the force here can also be relaxed  to any bounded force.
}

 {\remark
\lab{rmk2}

We now compare the conditions in the above theorem with the properties of the blow up solutions stated in \cite{Op}.

(1).   Per the construction, the order of the singularity is described  by the magnitude of negative powers of the variable $q=q(x_3, t)$ in the profile. Lower order singularities have less negative power.

 As described on p11-14 there,   the velocity $u$ in that solution is the sum of two functions: the axially symmetric background solution $u_B$ and the nonsymmetric perturbation $u^{ap}$ which is consisted of waves or  pulse functions called $w$ and a series of lower order corrections $u^{[j]}$:
\be
\lab{opai1}
u =\underbrace{ u^{co} + u^{sp}}_{u_B} + \underbrace{ w + \sum u^{[j]}}_{u^{ap}}.
\ee The superscript ``co'' stands for core, ``sp" stands for symmetric perturbation and ``ap'' for angled perturbation which is supported in annulus $\{X_a < X < X_b \}$.
Further, the back ground solution $u_B$ is the sum of the leading field  (core) $u^{co}=u^{co}(X, \eta)$   and an axially symmetric perturbation $u^{sp}$, which are constructed in Sections 3, 4 and Section 5 respectively. Note the $u^{co}$ only depend on $X$, $\eta$ and $u_B$ satisfies the profile \eqref{sscene} without the angle $\th$, c.f. (5.42) on p60. The wave function $w$ is constructed in Section 6 and it has a singularity of order $q^{-0.5 -0.5 h}$, c.f. p11 of \cite{Op}.
Comparing with the background solution $v_B$, this is $0.5h$ order smaller than the  singularity in the angular and axial velocity $v_\th$ and $v_3$ but $0.5 h$ bigger than that of the radial velocity $v_r$. The wave function is designed in such a way that $\nabla (w \otimes w)$ cancels the highest order term of the residue error from $u_B$. The singularity for $u^{[j]}$ is bounded from above by $q^{-K_0}$ on line 3 below (9.17), modulo a log term. We are not able to find an explicit expression for $K_0$. However, since $u^{[0]}$ is designed to cancel the residue error of $u_B+w$, whose singularity is already lower than that of the residue error of $u_B$, it is likely that $K_0$ is lower than $0.5 + 0.5 h$ which is the singularity order of $w$. The same reasoning can also be made for $u^{[j]}$. Therefore, it is likely that
the inequality $|E_2| \le C q^\delta$ in Condition (b) holds with $q^{-(0.5+h)} E_2$ representing the angular component of $u^{sp}+u^{ap}$. In fact, some descriptions in this regard are given on around (9.9) p106 which trace back to p85 around (7.35)-(7.36) in \cite{Op}. Note the space $W_\alpha$ defined in Definition 6.5 implies the functions there have an extra factor $\e^\alpha$ which is of the order $q^{h \alpha}$ by the convention in 6.1.

The other items of Conditions (a) and (b) are reasonable ones for profiles since singularities are mod out. According to Theorem 4.6 (i) and (iv), the leading profiles (for $u^{co}$) satisfy these items. By  (5.42) in Proposition 5.5, the back ground solution $u_B$ also  satisfy these items when $X$ is between two positive constants where $u^{sp}$ is supported. Therefore, $u_B$ satisfies these items.
It is possible that the profiles of the lower order terms $w+u^{ap}$ also satisfy these items. We have not finished verifying.

From this analysis, it seems possible that
the solution stated in \cite{Op} satisfies  Conditions (a) and (b) in Theorem \ref{th1.2}. However we can not draw a definite conclusion at this stage.

(2). The leading profile of the solution stated in \cite{Op} satisfies Condition (c).
See p33, (4.29) line 1 with $A=0.5+h$.  Since the perturbation happens in a region where the variable $X$ is bounded between two positive numbers, the whole solution also satisfies Condition (c).

(3). It is likely that the solution  in \cite{Op} satisfies Condition (d) for the same reason as  Conditions (a) and (b), explained in (1).

(4). In 3.3, p11 of \cite{Op}, it is described that the amplitude of the wave function $w$ in \eqref{opai1} is equivalent to $q^{-0.5-0.5h}$. This is the highest order perturbation added to the back ground solution $u_B$ whose radial component has singularity of order $q^{-0.5}$. Since the wave function have nontrivial radial components, see line -3, p11 and p71, it appears that
 the solution there satisfies Condition (e) with $\sig=0.5+0.5 h$. Note in the detailed construction in 6.1,  space time points are arranged in dyadic cubes such that: for a constant $C>1$, there exist $t_k \uparrow 0$, such that
\[
\limsup_{k}
\{ |R(\cdot, t_k)| \, : \, \, r <C|t_k|^{0.5}, \,  |x_3| < C |t_k|^{0.5-h} \}>0.
\]This just says that the scaled radial profile does not go to $0$ at singularity.
Therefore there exist $x_3(k)$, $\th_k$ and $r_k$ such that
\be
 \lab{limR>0}
 \lim_{k \to \infty} R(\frac{r^2_k}{ 2q(x_3(k), t_k)}, \eta(x_3(k), t_k), \th_k, t_k)>0.
 \ee
 By \eqref{taudenj0}, there is a constant $C$ such that
 \[
 \al
 &X_k \equiv \frac{r^2_k}{2 q(x_3(k), t_k)} \le C^2, \q |x_3(k)| \le C \tau_k^{0.5-h}, \\
 &\eta_k \equiv |\eta(x_3(k), t_k)| \le \frac{|x_3(k)|}{(\tau_k + |x_3(k)|^{1/(0.5-h)})^{0.5-h}} \le \frac{C}{(1+C^{1/(0.5-h)})^{0.5-h}}<1.
 \eal \]Hence there exist subsequences $X_k \to X_0<\infty$, $\eta_k \to \eta_0<1$ and $\th_k \to \th_0$. Taking limit in \eqref{limR>0}, we find
 \[
 R(X_0, \eta_0, \th_0, 0)>0.
 \]Therefore Condition (e) holds.

 Also it is mentioned, that there can be a log correction to the power of $q$, but this does not affect the proof of the theorem as explained in Remark \ref{rmk3.1} below.

One can also replace Condition (e) by the following one: $R(\cdot, t)$ is uniformly continuous for $(X, \eta, \th) \in [a, b] \times [-1, 1] \times [0, 2 \pi]$ and $v_r |_{t=0}$ does not vanish on open sets. Here $b>a>0$ are any finite numbers. The proof is almost identical.

(5). The solution stated in \cite{Op} satisfies condition \eqref{thrvth02}.
Indeed, by line -4 on p7 of that paper, the leading part of the angular velocity, denoted by $u^0_\th$ satisfies $u^0_\th = q^{-0.5-h} E$ with $A=0.5+h$  and $E$ being a positive function there. This and line 12 on p8 infers
\[
\al
\tau^h r u^0_\th(x, t) &= \tau^h r q^{-0.5-h} E(X, \eta)=  \tau^h r q^{-0.5-h} E(r^2/(2 q), \eta)\\
& =  (\tau/q)^h (r/q^{0.5}) E(r^2/(2 q), \eta)
\eal
\] Since $\tau = q$ for  $x_3=0$, the above implies
\[
\sup \tau^h r |u^0_\th|(\cdot, t) \ge \sup_{x_3=0} \tau^h r |u^0_\th|(\cdot, t)
=\sup_{X>0} \sqrt{2X} E(X, 0)>0,
\]by (3.6) of \cite{Op}. Hence \eqref{thrvth02} holds.

(6). All the boundedness conditions on functions and derivatives except (c) can be relaxed to local boundedness without affecting the proof.
}

\medskip

The idea of the proof of Theorem \ref{th1.1} is to use a blow up argument on the function $|t|^h \Gamma \equiv |t|^h r v_\th(x, t)$ which replaces the usual density function of the angular momentum $\Gamma=r v_\th$ that becomes unbounded as $t \uparrow 0$ under the stated scenario.  We will use the standard scaling for the Navier Stokes equations with spatial variables anchored at selected points on the $x_3$ axis, taking advantage of  the anisotropic nature of the velocity and the similar variables. In Theorem \ref{th1.1}, $v_r$ is critical but $v_\th$ and $v_3$ are super critical under the standard scaling. This causes the $x_3$ direction to collapse after scaling. It turns out the a suitable limit of $|t|^h \Gamma$ depends only on the $r$ and $t$ variable and it satisfies a one dimensional heat type equation which has an extra term involving the $\eta$ derivative of $E$ on one hand, but has no  pressure term on the other. When taking the limit, it is observed that all terms involving angular derivatives vanish eventually. We will argue that the extra term is $0$ at some special point which allows  the maximum principle to conclude the limit function is $0$. This is impossible since it arises as the angular profiles of the singularity. We notice that the forcing term does not play any role eventually since its scaling goes to $0$. For the proof of Theorem \ref{th1.2}, we argue that the divergence free condition coerces  the growth rate of $v_r$ to become critical. Then Theorem \ref{th1.1} will take over from here.

 We will use $C$ with or without index as positive constants that may change in value from line to line. If  variables need to be written out explicitly, we use $\f{\p \cdot}{\p \cdot}$  to denote partial derivatives. Sometimes, if there is no confusion, limits are taken for selected subsequences without explicitly mentioning so.

\section{ Proof of  theorem \ref{th1.1} }

The proof is divided into three steps.

{\it Step 1. introducing  a blow up procedure for $|t|^h r v_\th$ (change of variables by scaling and shifting)}.

\medskip

Suppose the singular solution satisfying the conditions of Theorem \ref{th1.1} exists so that \eqref{thrvth} holds, i.e.
\[
\limsup_{t \to 0^-} |t|^h \Vert r v_\th(\cdot, t) \Vert_\infty>0.
\]
Then, we can
pick a sequence of times $\{ t_k \}$, $t_k<0$ such that $t_k \uparrow 0$ as $k \to \infty$ and that
\be
\lab{thrvth2}
\limsup_{k  \to \infty} |t_k|^h \Vert r v_\th(\cdot, t_k) \Vert_\infty
=\limsup_{t \to 0^-} |t|^h \Vert r v_\th(\cdot, t) \Vert_\infty>0.
\ee
We will eventually reach a contradiction at the end of Step 3.

For $ t \le t_k$ and $x=(x_1, x_2, x_3) \in \mathbf{R}^3$, consider the functions
\be
\lab{defhk0}
H_k = \tau^h_k r v_\th(x, t), \quad r=\sqrt{x^2_1+x^2_2},
\ee and
\be
\lab{defjk0}
J_k = H_k + \e e^{-A \tau^{-0.5}_k r},
\ee where $\e>0$ is any small positive number and $A \ge 1$ is a large positive number to be specified later. We will prove, for large $k$,
\be
\lab{jk<2e0}
J_k \le 2 \e.
\ee

We  make the contrapositive assumption that, for a sequence of integers $k$ going to infinity,  $J_k$ reaches a maximum value greater than $2 \e$ at some point
\be
\lab{defp*}
p^*(k)=(r^*(k), \th(k), x^*_3(k), t^*(k))
\ee in the cylindrical coordinates. By renaming, we can take $k$ to be all positive integers. We next prove that $p^*(k)$ is localized around the origin in some sense.

From \eqref{defjk0}, at the point $p^*(k)$,
\be
\lab{defjk0.1}
2 \e \le J_k= (\tau_k/q)^{h} (r/q^{0.5}) \left( E( r^2/(2 q), \eta) +  E_2(r^2/(2 q), \eta, \theta, t)\right) + \e e^{-A \tau^{-0.5}_k r}.
\ee From \eqref{xetaq}, $q= \tau + x^2_3 q^{2 h} \ge x^2_3 q^{2h}$ so that $q \ge x^{2/(1-2h)}_3.$ This and Young's inequality yield, for a constant $C_1=C_1(h)$, that
\be
\lab{taudenj}
\tau + x^{2/(1-2h)}_3 \le q \le C_1 (\tau + x^{2/(1-2h)}_3)
\ee
Therefore
\be
\lab{tauk/q}
\f{\tau_k/C_1}{\tau +  x^{2/(1-2h)}_3} \le \f{\tau_k}{q} \le
\f{\tau_k}{\tau +  x^{2/(1-2h)}_3}.
\ee
From this and \eqref{defjk0.1}, we deduce, at $p^*_k$ that:
\be
\lab{e<jk}
\al
&\e \le J_k = (\tau_k/q)^{h} \underbrace{( r/q^{0.5})}_{\sqrt{X}}  \bigg( E(\underbrace{r^2 /(2 q)}_{X/2}, \cdot) +  E_2(\underbrace{r^2/(2q)}_{X/2}, \cdot)\bigg)\\
&\le
\left(\f{\tau_k}{\tau +  x^{2/(1-2h)}_3} \right)^h  \sup_X \sqrt{X} \, |(E+E_2)(X/2, \cdot)|.
\eal
\ee In the above and below "$\cdot$" stand for the variables $\eta, \th, t$.
By  Condition (c) in the theorem, there is a constant $C>0$ such that
\[
\sqrt{X} \, |(E+E_2)(X, \cdot)| \le C
\] for all $X>0$ uniformly. Therefore, \eqref{e<jk} infers $(q/\tau_k)^h \le C/\e$, and hence,
\be
\lab{t*jie0}
|t^*(k)| = \tau^*(k) \le (C/\e)^{1/h} \tau_k, \qq
|x^*_3(k)| \le (C/\e)^{(1-2h)/(2h)} \tau_k^{0.5-h}.
\ee

From $t \le t_k<0$ or equivalently $\tau \ge \tau_k>0$, one knows $\tau_k/q \le 1$ from
\eqref{tauk/q}. Thus \eqref{e<jk} implies
\[
\e \le( r/q^{0.5})   \, |(E+E_2)(r^2/(2 q), \cdot)|.
\]By Condition (c) again,
\[
\lim_{X \to \infty} \sqrt{X} \, |(E+E_2)(X, \cdot)|=0
\]uniformly. Therefore we know that there is constant $C_2(\e)>0$, independent of  $k$,  such that
\[
r/q^{0.5}
  \le C_2(\e).
\]By \eqref{taudenj} and \eqref{t*jie0}, we deduce
\be
\lab{r*jie0}
r^*(k) \le C_3(\e, h) \tau^{0.5}_k
\ee where $C_3(\e, h)$ is another constant, independent of  $k$.

The above argument shows that if the maximum of $J_k$ is greater than $2 \e$, it must occur at a point $p^*_k=(r^*(k), \th_k, x^*_3(k), t^*(k))$ which satisfy \eqref{t*jie0} and
\eqref{r*jie0}. The point is that $(r^*(k), \th_k, \tau_k^h x^*_3(k), t^*(k))$ is confined in a parabolic cube of size $C_4(\e, h) \tau^{0.5}_k$  where $C_4(\e, h)$ is independent of $k$.

In the cylindrical system, using the points $(r, \theta, x_3)=(0, 0, x^*_3(k))$ as the spatial anchors, we consider the following change of variables by shifting and scaling
\be
\lab{changvar}
\al
&t \Rightarrow t_k + |t_k| t=  t_k + \tau_k t,\\
&\tau \Rightarrow \tau_k + \tau_k \tau,\\
&r \Rightarrow \tau^{0.5}_k r, \qq x_3 \Rightarrow x^*_3(k) + \tau^{0.5}_k  x_3,\\
&q \Rightarrow q_k=\tau_k + \tau_k \tau + (x^*_3(k)+\tau^{0.5}_k x_3)^2 q^{2h}_k, \\
&X \Rightarrow X_k= \f{\tau_k r^2}{2 q_k}, \qq  \eta \Rightarrow \eta_k=\f{x^*_3(k) + \tau^{0.5}_k x_3}{q^{0.5-h}_k}.
\eal
\ee

  In the rest of the proof, $r, x_3, t, \tau$ will stand for the scaled variables (those after the arrows in \eqref{changvar}), unless stated otherwise. For simplicity of presentation, we will not give new names to the scaled variables.

Next we compute the limits of a few pertinent functions as $k \to \infty$.
For fixed $r, x_3, \tau$, by \eqref{taudenj} and \eqref{t*jie0},
\be
\lab{qkto0}
0 \le q_k  \le C_1 \left(\tau_k + \tau_k \tau + (x^*_3(k) +\tau^{0.5}_k x_3)^{2/(1-2h)} \right) = O(\tau_k)  \to 0,  \q k \to \infty.
\ee  Next
\be
\lab{tqxekjx0}
\al
 \f{\tau_k}{q_k} &= \f{\tau_k}{\tau_k + \tau_k \tau + (x^*_3(k)+\tau^{0.5}_k x_3)^2 q^{2h}_k } \\
&= \f{1}{1 +  \tau + (\tau_k/q_k)^{-2h} \left(\tau^{-0.5+h}_k x^*_3(k) + \tau^h_k x_3 \right)^2}.
\eal
\ee According to \eqref{t*jie0}, the sequence $\{\tau^{-0.5+h}_k x^*_3(k) \}$ is bounded. Hence we can find a subsequence, still indexed by $\{k\}$, and a finite number $x^\infty_3$ such that
\be
\lab{defx3wq}
\lim_{k \to \infty} \tau^{-0.5+h}_k x^*_3(k) = x^\infty_3.
\ee

Therefore, there exists a positive function $S=S(\tau)$ such that,
when $k \to \infty$,
\be
\lab{tqxekjx}
\al
& \f{\tau_k}{q_k} \to S(\tau)\\
& S(\tau)  = \f{1}{1 +  \tau + S^{-2h}(\tau) (x^\infty_3)^2},\\
&X_k = \f{\tau_k r^2}{2 q_k}  \to \f{r^2}{ 2} S(\tau),\\
& \eta_k = \f{x^*_3(k)+\tau^{0.5}_k x_3}{q^{0.5-h}_k} =
\left( \tau^{-0.5+h}_k x^*_3(k)+
  \tau^h_k x_3  \right) \left( \f{\tau_k}{q_k} \right)^{0.5-h} \to x^\infty_3 S^{0.5-h}(\tau).
  \eal
\ee It is clear that the above convergence can be regarded as in $C^{2, 1}_{x, t, loc}$ sense.

From line 2 of \eqref{tqxekjx0},
\[
\al
\f{\p}{\p x_3} \left( \f{\tau_k}{q_k} \right) & = \f{\p}{\p x_3} \left( 1 +  \tau + (\tau_k/q_k)^{-2h} \left(\tau^{-0.5+h}_k x^*_3(k) + \tau^h_k x_3 \right)^2 \right)^{-1}\\
&=- \left(\tau_k/q_k\right)^2 \,  \bigg[  (\tau_k / q_k )^{-2h} 2 \left(\tau^{-0.5+h}_k x^*_3(k)
 + \tau^h_k x_3 \right) \tau^h_k \\
 &\qq - 2h  (\tau_k / q_k )^{-2h-1} \f{\p}{\p x_3} \left( \f{\tau_k}{q_k} \right) \left(\tau^{-0.5+h}_k x^*_3(k) + \tau^h_k x_3 \right)^2 \bigg].
\eal
\]Hence
\be
\lab{p3t/q0}
\al
\f{\p}{\p x_3} \left( \f{\tau_k}{q_k} \right)& = - 2 \tau^{h}_k (\tau_k / q_k )^{2-2h}
\left(\tau^{-0.5+h}_k x^*_3(k) + \tau^h_k x_3 \right) \\
& \qquad \times \left[ 1- 2 h \left(\tau^{-0.5+h}_k x^*_3(k) + \tau^h_k x_3 \right)^2   (\tau_k / q_k )^{1-2h} \right]^{-1}\\
&=- 2 \tau^{h}_k \, (\tau_k / q_k )^{2-2h}
\left(\tau^{-0.5+h}_k x^*_3(k) + \tau^h_k x_3 \right) \\
& \qquad \times \left[ 1- 2 h \left(1-(1+\tau) (\tau_k/q_k) \right)   \right]^{-1}.
\eal
\ee In the last step, we just used \eqref{tqxekjx0}.

For fixed $\tau, x_3$ and large $k$, the denominator in the above equality is bounded away from $0$. Therefore  we deduce, from \eqref{p3t/q0} and \eqref{qkto0}
\be
\lab{p3t/q2}
\f{\p}{\p x_3} \left( \f{\tau_k}{q_k} \right) = O( \tau^{h}_k),
\q \f{\p}{\p x_3} q_k = O(\tau^{1+h}_k), \quad k \to \infty.
\ee
In addition
\be
\lab{p3t/q22}
\al
\lim_{k \to \infty}\tau^{-h}_k \f{\p}{\p x_3} \left( \f{\tau_k}{q_k} \right)
&=
- 2(S(\tau) )^{2-2h} x^\infty_3
 \left[ 1- 2 h + 2h (1+\tau) S(\tau) \right]^{-1} \\
&\equiv m_1(\tau).
\eal
\ee

 Consequently, by line 5 of \eqref{changvar}, for fixed $r>0$,
\be
\lab{p3Xk}
\f{\p}{\p x_3} X_k = \f{\p}{\p x_3} \left( \f{\tau_k}{q_k} \right) \f{r^2}{2} = O( \tau^{h}_k) r^2 = O( \tau^{h}_k), \quad k \to \infty
\ee and
\be
\lab{p3Xk2}
\al
\lim_{k \to \infty}\tau^{-h}_k \f{\p}{\p x_3} X_k &= \lim_{k \to \infty}\tau^{-h}_k \f{\p}{\p x_3} \left( \f{\tau_k}{q_k} \right) \f{r^2}{2}\\
&
 = - 2(S(\tau) )^{2-2h} x^\infty_3
 \left[ 1- 2 h + 2h (1+\tau) S(\tau) \right]^{-1} \f{r^2}{2}\\
 & \equiv m_2(\tau).
\eal
\ee

Finally, by \eqref{tqxekjx} line 4 and \eqref{p3t/q0},
\be
\lab{p3etak}
\al
&\f{\p}{\p x_3} \eta_k =
\tau^h_k  \left( \f{\tau_k}{q_k} \right)^{0.5-h}
+ \left( \tau^{-0.5+h}_k x^*_3(k)+
  \tau^h_k x_3  \right) (0.5-h) \left( \f{\tau_k}{q_k} \right)^{-0.5-h} \f{\p}{\p x_3} (\tau_k/q_k)\\
&= \tau^h_k  \left( \f{\tau_k}{q_k} \right)^{0.5-h} - 2 \tau^{h}_k (\tau_k / q_k )^{1.5-3h}
\left(\tau^{-0.5+h}_k x^*_3(k) + \tau^h_k x_3 \right) \\
& \qquad \times \left[ 1- 2 h \left(1-(1+\tau) (\tau_k/q_k) \right)   \right]^{-1}
\left( \tau^{-0.5+h}_k x^*_3(k)+
  \tau^h_k x_3  \right) (0.5-h).
\eal
\ee Therefore,
\be
\lab{p3etak2}
\al
&\f{\p}{\p x_3} \eta_k
=O( \tau^{h}_k), \quad k \to \infty;\\
&\lim_{k \to \infty}\tau^{-h}_k \f{\p}{\p x_3} \eta_k \\
&
\qq =S(\tau)^{0.5-h} - 2 S(\tau)^{1.5-2h} x^\infty_3
\left[ 1- 2 h + 2h (1+\tau) S(\tau)   \right]^{-1}
x^\infty_3 (0.5-h) \\
&\qq \equiv m_3(\tau).
\eal
\ee Again, the above convergence can be regarded as in $C^{2, 1}_{x, t, loc}$ sense.

\medskip

{\it Step 2.  Scaled equation for $r v_\th \tau^{h}_k$.}

\medskip

Using line 2 of \eqref{eqasns}, in the original variables, the density of angular momentum $r v_\th$ satisfies the equation
\be
\lab{eqrvth}
\al
&\p^2_r (r v_\th) - \f{1}{r} \p_r( r v_\th) + \p^2_{x_3} (r v_\th) + \f{1}{r^2} \p^2_\th(r v_\th)\\
&\quad -v_r \p_r (r v_{\theta})  - \f{1}{r} v_\th \p_\theta (r v_{\theta})- v_3 \p_{x_3}  (r v_{\theta}) + \f{2}{r} \p_\theta  v_r -  \p_\theta P-
\partial_t (r v_{\theta})= r f_\th
\eal
\ee We will multiply the above equation by $\tau_k^{1+h}$ and convert to the scaled variables in \eqref{changvar}. There are 11 terms to deal with.

Using the given profile for $v_\th$, we convert
\[
\al
&T_{1k} \equiv \tau_k^{1+h} \p^2_r (r v_\th)\\
&\Rightarrow T_{1k} = \tau_k^{1+h} \f{\p^2 \left(\tau^{0.5}_k r \f{1}{q^{0.5+h}_k} \left( E(X_k, \eta_k) +  E_2(X_k, \eta_k, \theta, s_k)\right) \right)}{ \p(\tau^{0.5}_k r)^2}\\
&=\f{\p^2 \left( r (\tau_k/q_k)^{0.5+h} \left( E(X_k, \eta_k) +  E_2(X_k, \eta_k, \theta, s_k)\right) \right)}{ \p r^2}.
\eal
\]Here and later in the proof,
\be
\lab{defsk}
s_k=t_k+\tau_k t = t_k - \tau_k \tau
\ee with $t$ being the new time variable. Let us write the function $H_k$ in terms of the scaled variables $r, x_3, \tau$ and $\theta$:
\be
\lab{defhk}
\al
H_k&=r (\tau_k/q_k)^{0.5+h} \left( E(X_k, \eta_k) +  E_2(X_k, \eta_k, \theta, s_k)\right)\\
&=  (\tau_k/q_k)^{0.5+h} r \left( E\left(\f{\tau_k}{q_k} \f{r^2}{2}, (\tau^{-0.5+h}_k x^*_3(k)+\tau^h_k x_3) \left( \f{\tau_k}{q_k} \right)^{0.5-h}\right) +  E_2(X_k, \eta_k, \theta, s_k)\right).
\eal
\ee The last two identities infer
\be
\lab{t1scal}
T_{1k} \equiv \tau_k^{1+h} \p^2_r (r v_\th) \Rightarrow T_{1k}=\f{\p^2}{\p r^2} H_k.
\ee Using \eqref{tqxekjx} and Condition (b) of the theorem, we have, in the $C^0_{X, \eta, \th, loc}$ sense,
\be
\lab{hktoh}
\lim_{k \to \infty} H_k = H,
\ee where
\be
\lab{defH}
H \equiv S(\tau)^{0.5+h} r E(\f{r^2}{2} S(\tau), x^\infty_3 S(\tau)^{0.5-h}).
\ee Notice that $H$ is independent of the angle $\th$ and $x_3$.

Next
\[
\al
T_{2k} \equiv \tau_k^{1+h} \f{1}{r} \p_r( r v_\th) &\Rightarrow
\tau_k^{1+h} \f{1}{ \tau_k^{0.5} r}
\f{\p \left(\tau^{0.5}_k r \f{1}{q^{0.5+h}_k} \left( E(X_k, \eta_k) +  E_2\right) \right)}{ \p(\tau^{0.5}_k r)}\\
&=\f{1}{r}
\f{\p \left( r \left(\tau_k/q_k \right)^{0.5+h} \left( E(X_k, \eta_k) +  E_2\right) \right)}{ \p r} = \f{1}{r} \f{\p}{\p r} H_k.
\eal
\]i.e.
\be
\lab{t2scal}
T_{2k} \equiv \tau_k^{1+h} \f{1}{r} \p_r( r v_\th)  \Rightarrow T_{2k} = \f{1}{r} \f{\p}{\p r} H_k.
\ee

Now we deal with
\be
\lab{t3scal}
\al
&T_{3k} \equiv \tau_k^{1+h} \p^2_{x_3} (r v_\th) \\
&\Rightarrow
T_{3k} =\tau_k^{1+h}  \f{\p^2 \left(\tau^{0.5}_k r \f{1}{q^{0.5+h}_k} \left( E(X_k, \eta_k) +  E_2\right) \right)}{ \p(\tau^{0.5}_k x_3)^2}\\
&=\f{\p^2 \left( r (\tau_k/q_k)^{0.5+h} \left( E(X_k, \eta_k) +  E_2(X_k, \eta_k, \theta, s_k)\right) \right)}{ \p x_3^2}\\
&=\f{\p^2}{\p x_3^2} H_k.
\eal
\ee By differentiating each term in line 2 of \eqref{t3scal}, one can show that
$\f{\p^2}{\p x_3^2} H_k \to 0$ uniformly on compact sets. Since the calculation is tedious, in the next step, we will instead use an integral argument to show this term vanishes as $k \to \infty$.

Now we treat the 4th term on the left side of the equation for $r v_\th$.
\be
\lab{t4scal}
\al
&T_{4k} \equiv \tau_k^{1+h} \f{1}{r^2} \p^2_\th(r v_\th)\\
&\Rightarrow T_{4k} =\tau_k^{1+h} \f{1}{\tau_k r^2} \f{\p^2 \left(\tau^{0.5}_k r  \left(\f{1}{q^{0.5+h}_k} \left( E(X_k, \eta_k) +  E_2\right) \right)\right)}{ \p \th^2}\\
&=\f{1}{r^2}  \f{\p^2 \left( r \left(\tau_k/q_k \right)^{0.5+h} \left( E(X_k, \eta_k) +  E_2\right) \right)}{ \p \th^2}\\
&=\f{1}{r^2}  \f{\p^2 H_k}{ \p \th^2}
\eal
\ee

Next comes the radial drift term.
\be
\lab{t5scal}
\al
&T_{5k} \equiv \tau_k^{1+h} v_r \p_r( r v_\th) \\
&\Rightarrow
T_{5k} = \tau_k^{1+h} \f{1}{ q_k^{0.5}}  R(X_k, \eta_k, \th, s_k)
\f{\p \left(\tau^{0.5}_k r \f{1}{q^{0.5+h}_k} \left( E(X_k, \eta_k) +  E_2\right) \right)}{ \p(\tau^{0.5}_k r)}\\
&=(\tau_k/q_k)^{0.5}  R(X_k, \eta_k, \th, s_k)
\f{\p \left(  \left(\tau_k/q_k \right)^{0.5+h} r \left( E(X_k, \eta_k) +  E_2\right) \right)}{ \p r}\\
& =(\tau_k/q_k)^{0.5} R(X_k, \eta_k, \th, s_k)  \f{\p}{\p r} H_k\\
&\to S(\tau)^{0.5+h} R(\f{r^2}{2} S(\tau), x^\infty_3 S(\tau)^{0.5-h}, \th, 0) \f{\p}{\p r}  H, \q k \to \infty.
\eal
\ee Here the function $H$ is defined in \eqref{defH} and we have used Conditions (a) and (b) of the theorem. The above convergence is in the sense of uniform convergence in compact sets. Note, by \eqref{defsk},  $s_k=t_k+\tau_k t \to 0$ as $k \to \infty$ for any fixed scaled time $t$.

For the angular drift term,
\be
\lab{t6scal}
\al
&T_{6k} \equiv \tau_k^{1+h} v_\th \p_\th v_\th \\
&\Rightarrow
T_{6k}=\tau_k^{1+h} \f{1}{ q_k^{0.5+h}} \left(E(X_k, \eta_k) +  E_2 \right)
\f{\p \left(\f{1}{q^{0.5+h}_k} \left( E(X_k, \eta_k) +  E_2\right) \right)}{ \p \th}\\
&= \tau^{-h}_k \left(\tau_k/q_k \right)^{0.5+h} \left(E(X_k, \eta_k) +  E_2\right)
\f{\p \left(  \left(\tau_k/q_k \right)^{0.5+h} \left( E(X_k, \eta_k) +  E_2\right) \right)}{ \p \th}\\
& = \tau^{-h}_k r^{-2} H_k \f{\p}{\p \th} H_k.
\eal
\ee
 Due to the factor $\tau^{-h}_k$,  this term does not converge as $k \to \infty$. However its integral with a compactly supported axially symmetric function is $0$. This property will make this term vanishes after taking the limit. Note this is no singularity at $r=0$ since $E=0$ when $r=0$.

Our next target is the axial drift term which needs special attention.
\be
\lab{t7scal}
\al
&T_{7k} \equiv \tau_k^{1+h} v_3 \p_{x_3} (r v_\th)\\
&\Rightarrow
T_{7k} =\tau_k^{1+h} \f{1}{q^{0.5+h}_k} U(X_k, \eta_k, \th, s_k) \f{\p \left(\tau^{0.5}_k r \f{1}{q^{0.5+h}_k} \left( E(X_k, \eta_k) +  E_2\right) \right)}{ \p(\tau^{0.5}_k x_3)}\\
&=\tau^{-h}_k  (\tau_k/q_k)^{0.5+h} U(X_k, \eta_k, \th, s_k) r \underbrace{\f{\p \left(  (\tau_k/q_k)^{0.5+h} ( E(X_k, \eta_k) + E_2)  \right)}{ \p x_3}}_{I_1}.
\eal
\ee Here and in the next few lines, we use $\cdot$, $\cdot$ to denote $\theta$ and $s_k$ when no confusion arises.  We need to cancel the diverging term $\tau^{-h}_k$ in the preceding line.
Note that
\[
\al
&I_1=\f{\p \left(  (\tau_k/q_k)^{0.5+h} \left( E(X_k, \eta_k) +  E_2(X_k, \eta_k, \th, s_k)\right) \right)}{ \p x_3}\\
&=(0.5+h) (\tau_k/q_k)^{-0.5+h} \f{\p}{\p x_3} (\tau_k/q_k) \left( E(X_k, \eta_k) +  E_2(X_k, \eta_k, \th, s_k)\right)\\
&\q + (\tau_k/q_k)^{0.5+h} \left( \f{\p}{\p X} E( X_k, \eta_k) \f{\p X_k}{\p x_3} +
\f{\p}{\p \eta} E( X_k, \eta_k) \f{\p \eta_k}{\p x_3} \right)\\
&\q + (\tau_k/q_k)^{0.5+h} \left( \f{\p}{\p X} E_2( X_k, \eta_k, \th, s_k) \f{\p X_k}{\p x_3} +
 \f{\p}{\p \eta} E_2 \f{\p \eta_k}{\p x_3} \right).
\eal
\]By \eqref{qkto0} and Condition (b), we have
\be
\lab{e2to0}
|E_2( X_k, \eta_k, \th, s_k)| \le C q^\delta_k = O(\tau^\delta_k),
\ee

Using \eqref{p3t/q2}, \eqref{p3Xk} and \eqref{p3etak2}, we deduce
Thus
\be
\lab{i12}
\al
&I_1
=(0.5+h) (\tau_k/q_k)^{-0.5+h} \f{\p}{\p x_3} (\tau_k/q_k)  E(X_k, \eta_k) \\
&\q + (\tau_k/q_k)^{0.5+h} \left( \f{\p}{\p X} (E+E_2) \f{\p X_k}{\p x_3} +
\f{\p}{\p \eta} (E+E_2) \f{\p \eta_k}{\p x_3} \right) + o(1) \tau^h_k.
\eal
\ee Here and later in the proof, $o(1)$ stands for terms that converge to $0$ uniformly on compact sets as $k \to \infty$.
Substituting this into \eqref{t7scal}, we can cancel the term $\tau^{-h}_k$ to deduce
\be
\lab{t7scal2}
\al
&T_{7k} \equiv \tau_k^{1+h} v_3 \p_{x_3} (r v_\th)\\
&\Rightarrow T_{7k} =
  (\tau_k/q_k) U(X_k, \eta_k,  \th, s_k) r \bigg[
 (0.5+h) (\tau_k/q_k)^{-0.5+h} \tau_k^{-h} \f{\p}{\p x_3} (\tau_k/q_k)  E(X_k, \eta_k) \\
&\q + (\tau_k/q_k)^{0.5+h} \left( \f{\p}{\p X} (E+E_2) \tau_k^{-h} \f{\p X_k}{\p x_3} +
\f{\p}{\p \eta} (E+E_2) \tau_k^{-h} \f{\p \eta_k}{\p x_3} \right) + o(1) \bigg].
\eal
\ee This, \eqref{p3t/q22}, \eqref{p3Xk2} and \eqref{p3etak2} imply the uniform  convergence on compact sets:
\be
\lab{t7kjx}
\al
&\lim_{k \to \infty} T_{7k}\\
& =
  S(\tau) U r \bigg[
 (0.5+h) S(\tau)^{-0.5+h} m_1(\tau)  E \\
 &\qq \qq + S(\tau)^{0.5+h} \left( \f{\p (E+E_2)}{\p X}  m_2(\tau) +
\f{\p (E+E_2)}{\p \eta}  m_3(\tau) \right)\bigg].
\eal
\ee In the above the arguments for the functions $U$ and $E+E_2$ are $(\f{r^2}{2} S(\tau), x^\infty_3 S(\tau)^{0.5-h}, \th, 0)$ again.

Next we deal with the 8th term.
\be
\lab{t8scal}
\al
&T_{8k} \equiv \tau_k^{1+h} \f{2}{r} \p_\th v_r \\
&\Rightarrow T_{8k} = \tau_k^{1+h} \f{1}{\tau^{0.5}_k r} \frac{1}{q^{0.5}_k} \f{\p}{\p \th} R(X_k, \eta_k, \th, s_k)\\
& = \tau_k^{h} \f{2}{ r} (\tau_k/q_k)^{0.5} \f{\p}{\p \th} R(X_k, \eta_k, \th, s_k) \to 0, \quad k \to \infty.
\eal
\ee The above convergence is  uniform due to the factor $\tau_k^h$.

Now we come to the pressure term. Although we do not have any estimate on the pressure, the pressure term contains the $\th$ derivative. When testing against compactly supported axially symmetric functions, this  makes it vanishing after taking the limit. So we simply write
\be
\lab{t9scal}
T_{9k} \equiv \tau_k^{1+h} \p_\th P \Rightarrow T_{9k}=\tau_k^{1+h} \p_\th P(X_k, \eta_k, \th, s_k).
\ee For the same reason, we even do not need uniform boundedness of $\p R/\p \th$ since the same type of integrals testing of $T_{8k}$ is $0$.

The time derivative term is treated as follows
\be
\lab{t10scal}
\al
&T_{10k} \equiv \tau_k^{1+h} \p_t( r v_\th)\\
&\Rightarrow T_{10k}
=\tau_k^{1+h} \f{\p}{\p (t_k + \tau_k t)} \left(\tau^{0.5}_k r  \f{1}{q^{0.5+h}_k} (E+E_2)\right)\\
&= \f{\p}{\p t} \left(   (\tau_k/q_k)^{0.5+h} r (E+E_2)(X_k, \eta_k, \th, s_k)\right) = \f{\p}{\p t} H_k.
\eal
\ee  Finally since the forcing term is bounded, $\tau^{1+h}_k f_\th$ converges to $0$ uniformly as $k \to \infty$.

Gathering the terms $T_{ik}$, $i=1, ..., 10$ and the forcing term, we obtain the following equation for $H_k=H_k(r, x_3, \th, \tau)$.
\be
\lab{eqhk}
\al
&\underbrace{\f{\p^2}{\p r^2} H_k}_{T_{1k}}- \underbrace{\f{1}{r} \f{\p}{\p r} H_k}_{T_{2k}} + \underbrace{\f{\p^2}{\p x_3^2} H_k}_{T_{3k}} + \underbrace{\f{1}{r^2}  \f{\p^2}{ \p \th^2} H_k}_{T_{4k}} - \underbrace{(\tau_k/q_k)^{0.5} R(X_k, \eta_k, \th, s_k)  \f{\p}{\p r} H_k}_{T_{5k}}\\
&\q - \underbrace{\tau^{-h}_k r^{-2} H_k \f{\p}{\p \th} H_k}_{T_{6k}}
- \q \underbrace{ T_{7k}}_{\text{given below}} \q
 + \underbrace{\tau_k^{h} \f{2}{ r} (\tau_k/q_k)^{0.5} \f{\p}{\p \th} R(X_k, \eta_k, \th, s_k)}_{T_{8k}}\\
 &\q
-\underbrace{\tau_k^{1+h} \p_\th P(X_k, \eta_k, \th, s_k) }_{T_{9k}}
- \underbrace{\f{\p}{\p t} H_k}_{T_{10k}} = \tau^{1.5+h}_k r f_\th,
\eal
\ee where, as given in \eqref{t7scal2},
\[
\al
T_{7k} &=
  (\tau_k/q_k) U(X_k, \eta_k, \th, s_k)) r \bigg[
 (0.5+h) (\tau_k/q_k)^{-0.5+h} \tau_k^{-h} \f{\p}{\p x_3} (\tau_k/q_k)  E(X_k, \eta_k) \\
&\q + (\tau_k/q_k)^{0.5+h} \left( \f{\p}{\p X} (E+E_2) \tau_k^{-h} \f{\p X_k}{\p x_3} +
\f{\p}{\p \eta} (E+E_2) \tau_k^{-h} \f{\p \eta_k}{\p x_3} \right) + o(1) \bigg].
\eal
\]In the above equation, we move most terms which will vanish eventually to the right hand side, obtaining
\be
\lab{eqhk2}
\al
L_k  H_k &\equiv \f{\p^2}{\p r^2} H_k - \f{1}{r} \f{\p}{\p r} H_k - \f{\p}{\p t} H_k -
(\tau_k/q_k)^{0.5} R(X_k, \eta_k, \th, s_k)  \f{\p}{\p r} H_k\\
&\qquad  - (\tau_k/q_k) U(X_k, \eta_k, \th, s_k)  r \bigg[
 (0.5+h) (\tau_k/q_k)^{-0.5+h} \tau_k^{-h} \f{\p}{\p x_3} (\tau_k/q_k)  E(X_k, \eta_k) \\
&\qq \qq + (\tau_k/q_k)^{0.5+h} \left( \f{\p}{\p X} (E+E_2) \tau_k^{-h} \f{\p X_k}{\p x_3} +
\f{\p}{\p \eta}(E+E_2) \tau_k^{-h} \f{\p \eta_k}{\p x_3} \right) + o(1) \bigg]\\
&= - \f{\p^2}{\p x_3^2} H_k - \f{1}{r^2}  \f{\p^2}{ \p \th^2} H_k
+ \tau^{-h}_k r^{-2} H_k \f{\p}{\p \th} H_k - \tau_k^{h} \f{2}{ r} (\tau_k/q_k)^{0.5} \f{\p}{\p \th} R(X_k, \eta_k, \th, s_k)\\
&\qq + \tau_k^{1+h} \p_\th P(X_k, \eta_k, \th, s_k) + \tau^{1.5+h}_k r f_\th.
\eal
\ee

\medskip

{\it Step 3.  taking limits.}

\medskip

Recall from \eqref{defjk0}  and \eqref{defhk},  in the scaled variables the function $J_k$ becomes
\be
\lab{defjk}
\al
&J_k = H_k + \e e^{-A  r}=\\
&\left(\f{\tau_k}{q_k}\right)^{0.5+h} r \left( E\left(\f{\tau_k}{q_k} \f{r^2}{2}, (\tau^{-0.5+h}_k x^*_3(k)+\tau^h_k x_3) \left( \f{\tau_k}{q_k} \right)^{0.5-h}\right) +  E_2(X_k, \eta_k, \theta, s_k)\right) + \e e^{-A  r}.
\eal
\ee
 Direct computation from \eqref{eqhk2} shows
\be
\lab{eqjk}
\al
&L_k  J_k \equiv \f{\p^2}{\p r^2} J_k - \f{1}{r} \f{\p}{\p r} J_k - \f{\p}{\p t} J_k -
(\tau_k/q_k)^{0.5} R(X_k, \eta_k, \th, s_k) \f{\p}{\p r} J_k\\
&\qquad   - (\tau_k/q_k) U(X_k, \eta_k, \th, s_k) r \bigg[
 (0.5+h) (\tau_k/q_k)^{-0.5+h} \tau_k^{-h} \f{\p}{\p x_3} (\tau_k/q_k)  E(X_k, \eta_k) \\
&\qq \qq + (\tau_k/q_k)^{0.5+h} \left( \f{\p}{\p X} (E+E_2) \tau_k^{-h} \f{\p X_k}{\p x_3} +
\f{\p}{\p \eta} (E+E_2) \tau_k^{-h} \f{\p \eta_k}{\p x_3} \right) + o(1) \bigg]
\\
&= - \f{\p^2}{\p x_3^2} H_k - \f{1}{r^2}  \f{\p^2}{ \p \th^2} H_k
+ \tau^{-h}_k r^{-2} H_k \f{\p}{\p \th} H_k - \tau_k^{h} \f{2}{ r} (\tau_k/q_k)^{0.5} \f{\p}{\p \th} R(X_k, \eta_k, \th, s_k)\\
&\qq + \tau_k^{1+h} \p_\th P(X_k, \eta_k, \th, s_k) + \tau^{1.5+h}_k  r f_\th\\
&\qq + \e A e^{-A r} \left( A +\f{1}{r} + R(X_k, \eta_k, \th, s_k) \right).
\eal
\ee

Let $\phi=\phi(r, x_3, t)$ be a smooth, compactly supported in $\{ r>0 \}$, axially symmetric function which vanishes for $t \le -T$ where $T>0$ is a given positive number. Testing \eqref{eqjk} with $\phi$ yields
\be
\lab{eqjkint}
\al
&\int^0_{-T} \int^{2\pi}_0 \int \int \phi \, (L_k J_k) r dr dx_3 d\th dt \\
&=\int^0_{-T} \int^{2\pi}_0 \int \int \phi \, \e A e^{-A r} \left( A +\f{1}{r} +
 R(X_k, \eta_k, \th, s_k)\right) dr dx_3 d\th dt\\
&\qq -\int^0_{-T}  \int^{2\pi}_0 \int \int \phi \, \f{\p^2}{\p x_3^2} H_k r dr dx_3 d\th dt
+ \tau^{1.5+h}_k \int^0_{-T}  \int^{2\pi}_0 \int \int \phi \, f_\th r^2 dr dx_3 d\th dt.
\eal
\ee In reaching this identity, we have used crucially the fact that all terms involving $\f{\p}{\p \th}$ on the righthand side of \eqref{eqjk} have $0$ integrals. For example
\[
\int^{2\pi}_0 \int \int \phi(r, x_3)  H_k \f{\p}{\p \th} H_k r dr dx_3 d\th =0.
\]

We need to treat the time derivative term for $J_k$ inside $L_k J_k$ more carefully since we are not assuming $\f{\p}{\p t} H_k$ converges point-wise. From \eqref{defjk} and
\eqref{defhk},
\be
\lab{ptjk}
\f{\p}{\p t} J_k = \f{\p}{\p t} H_k=  r \f{\p}{\p t}\left[(\tau_k/q_k)^{0.5+h} \left( E(X_k, \eta_k) +  E_2(X_k, \eta_k, \theta, s_k)\right)\right].
\ee Therefore, writing $dx = r dr dx_3 d \th$ for simplicity, we find
\be
\lab{iptjk}
\al
&I_{\p_t J_k} \equiv \int^0_{-T} \int^{2\pi}_0 \int \int \phi \, \f{\p}{\p t} J_k r dr dx_3 d\th dt\\
&=\int^0_{-T} \int \phi \,   r \f{\p}{\p t}[(\tau_k/q_k)^{0.5+h}  E(X_k, \eta_k)]  dx dt \\
&\qq + \int^0_{-T} \int \phi \,   r \f{\p}{\p t}[(\tau_k/q_k)^{0.5+h}  E_2] dx dt\\
&=-\int^0_{-T} \int \f{\p}{\p t} \phi \,  \underbrace{ r [(\tau_k/q_k)^{0.5+h}  E(X_k, \eta_k)]}_{H_k}  d  x dt + \int \phi \,   r [(\tau_k/q_k)^{0.5+h}  E(X_k, \eta_k)] \big |_{t=0} dx \\
&\qq - \int^0_{-T} \int \f{\p}{\p t} \phi \,   r [(\tau_k/q_k)^{0.5+h}  E_2]  dx dt + \int \phi \,   r [(\tau_k/q_k)^{0.5+h}  E_2 ] \big |_{t=0} dx.
\eal
\ee   Recall by Condition (b) of the theorem, \eqref{e2to0} holds, and hence $E_2$ goes to $0$ in $C^0_{X, \eta, \th, loc}$ sense as $k \to \infty$.

Since $\phi$ is compactly supported  and by \eqref{hktoh}, as $k \to \infty$,
\[
H_k = H_k(r, x_3, \th, t) \to H=H(r, \tau)=S(\tau)^{0.5+h} r E(\f{r^2}{2} S(\tau), x^\infty_3 S(\tau)^{0.5-h}).
\] in $C^{0}_{loc}$ sense, we deduce
\be
\lab{iptjk2}
\al
\lim_{k \to \infty} I_{\p_t J_k} &= -\int^0_{-T} \int \f{\p}{\p t} \phi \,  H
  d x dt + \int \phi \,   H \big |_{t=0} dx \\
&=-\int^0_{-T} \int \f{\p}{\p t} \phi \, ( H + \e e^{-Ar})  dx dt + \int \phi \,   ( H + \e e^{-Ar}) \big |_{t=0} dx\\
&=\int^0_{-T} \int^{2\pi}_0 \int \int  \phi \,  \f{\p}{\p t}  J r dr dx_3 d\th dt.
\eal
\ee In the above we have written
\be
\lab{defJ}
J=J(r, t) \equiv H + \e e^{-A r}= S(\tau)^{0.5+h} r E(\f{r^2}{2} S(\tau), x^\infty_3 S(\tau)^{0.5-h}) + \e e^{-A r},
\ee and we only need the convergence of $H_k$ is $C^0_{loc}$ sense.

Similarly, we obtain
\be
\al
&\lim_{k \to \infty} \int^{2\pi}_0 \int \int \phi \, \f{\p^2}{\p x_3^2} H_k  r dr dx_3 d\th = \lim_{k \to \infty} \int^{2\pi}_0 \int \int \left(\f{\p^2}{\p x_3^2} \phi \right)\,  H_k  r dr dx_3 d\th\\
&= \int^{2\pi}_0 \int \int \left( \f{\p^2}{\p x_3^2} \phi \right)\,  H(r, \tau) r dr dx_3 d\th =0.
\eal
\ee

Under Conditions (a) and (b) of the theorem, we can now take the limit in \eqref{eqjkint}, using the limits of $T_{ik}$, $i=1, 2, ...10$,  to reach
\be
\lab{eqjjifen}
\al
\int^{0}_{-T}\int \phi(r, x_3, t) &\bigg(
\f{\p^2}{\p r^2} J - \f{1}{r} \f{\p}{\p r} J - \f{\p}{\p t} J -
S(\tau)^{0.5} R \f{\p}{\p r} J\\
&- S(\tau) U r \bigg[
 (0.5+h) S(\tau)^{-0.5+h} m_1(\tau)  E \\
 &\qq + S(\tau)^{0.5+h} \big( \underbrace{\f{\p E}{\p X} \, m_2(\tau)}_{Y_1} +
\f{\p E}{\p \eta}  m_3(\tau) \big)\bigg]
  \bigg)
dx dt\\
&=\int^{0}_{-T} \int  \phi(r, x_3, t) \, \e A e^{-A r} \left( A +\f{1}{r} + R \right)  dx dt.
\eal
\ee In the above, the limits in \eqref{tqxekjx} were used repeatedly and all terms involving $E_2$ vanishes. Note, there is no need for information on the  spatial second derivative derivatives for $E, E_2$. For instance, by Condition (b) of the theorem (see also \eqref{e2to0}),
\[
\int \phi \Delta E_2 dx= \int E_2 \Delta \phi dx \to 0, \q k \to \infty,
\]and, by Condition (a) and \eqref{defjk},
\[
\al
&\lim_{k \to \infty} \int^{0}_{-T}\int \phi \Delta J_k dxdt=
\lim_{k \to \infty} \int^{0}_{-T}\int \Delta \phi  J_k dxdt\\
&=\int^{0}_{-T}\int \Delta \phi  J dxdt.
\eal
\]

The convergence to  the term $Y_1$ also needs a little explanation since we only assumed $\sqrt{X} \p_X E$ is bounded in Condition (a). By \eqref{p3Xk} and \eqref{p3Xk2}, the terms $\p_{x_3} X_k$ and $m_2(\tau)$ all contain the factor $r^2$. Therefore, the following inequality holds (from 3r line of \eqref{eqjk})
\[
\left| \f{\p}{\p X} E \tau_k^{-h} \f{\p X_k}{\p x_3} \right| \le C X_k^{-1/2} r^2
=C (q_k/\tau_k)^{0.5} r^{-1} r^2 = C (q_k/\tau_k)^{0.5} r.
\]So the left hand side is locally bounded. The dominated convergence theorem infers its convergence  to
$
| \p_X (E+E_2) | m_2(\tau)
$ in the integral sense. This observation will be used later around \eqref{detaEe2}.

Since the integrands in \eqref{eqjjifen} except $R$ and $U$ are independent of $\th$ and only $\phi$ depends on $x_3$, we see that $J=J(r, t)$ is a bounded, ancient solution to the following equation: for $r>0$,
\be
\lab{eqjgudai}
\al
&\f{\p^2}{\p r^2} J - \f{1}{r} \f{\p}{\p r} J - \f{\p}{\p t} J -
S(\tau)^{0.5} \widetilde{R} \, \f{\p}{\p r} J\\
&   - \underbrace{ S(\tau)\widetilde{U} r \big[
 (0.5+h) S(\tau)^{-0.5+h} m_1(\tau)  E + S(\tau)^{0.5+h} \big( \f{\p E}{\p X}  \, m_2(\tau) +
\f{\p E}{\p \eta}  \, m_3(\tau) \big)\big] }_{I_2}
\\
&= \e A e^{-A r} \left( A +\f{1}{r} + \widetilde{R} \right).
\eal
\ee In the above
\[
\al
\widetilde{R} &= (1/2 \pi) \int^{2\pi}_0 R(\cdot, \cdot, \th, 0) d\th, \q \widetilde{U} = (1/2 \pi) \int^{2\pi}_0 U(\cdot, \cdot, \th, 0) d\th,
 \\&(\cdot, \cdot)=(\f{r^2}{2} S(\tau), x^\infty_3 S(\tau)^{0.5-h}).
\eal
\]
Since $R$ is a bounded function, we can and do choose the constant $A$ large enough so that the right hand side is positive a very point $(r, t) \in (0, \infty) \times (-\infty, 0]$.

If the under braced term $I_2$ were not there, we could then use the maximum principle to show that $J \le 0$. To get around this difficulty, let us go back to the function $J_k$ and traces where its maximum can occur.

Recall from Step 1 that at $p^*_k=(r^*(k), \th_k, x^*_3(k), t^*(k))$,   the maximum point of $J_k$  in the original variables, inequality \eqref{t*jie0} and
\eqref{r*jie0} hold. This means that
\[
(r^*(k), \th_k, \tau_k^h x^*_3(k), t^*(k))
\]is confined in a parabolic cube of size $C_4(\e, h) \tau^{0.5}_k$  where $C_4(\e, h)$ is independent of $k$.
Now we use the scaled variables $(r, \th, x_3, t)$ in \eqref{changvar}. Then
\be
\lab{p*2}
p^*_k = (r, \th, x_3, t)=(\tau^{-0.5}_k r^*(k), \th_k, 0, \tau^{-1}_k t^*(k)-1).
\ee Hence there exists a subsequence, still denoted by $p^*_k$ and a number $r_0 > 0$, such that
\be
\lab{p*2jx}
p_\infty \equiv \lim_{ k \to \infty} p^*_k  =(r_0, 0, 0, 0).
\ee Here we have written the limits of $\th_k$ and $\tau^{-1}_k t^*(k)-1$ as $0$ for simplicity.
Note that $r_0$ is bounded away from $0$ since the maximum value of $J_k$, which is greater than $2 \e$,  would converge to $\e$.

From \eqref{defhk} and \eqref{defjk}, in the scaled variables,
\be
\lab{defjk2}
J_k=r (\tau_k/q_k)^{0.5+h} \left( E(X_k, \eta_k) +  E_2(X_k, \eta_k, \theta, s_k)\right) + \e e^{-A r}
\ee where $X_k= r^2 \tau_k/(2 q_k)$ and $\eta_k=(\tau^{-0.5+h}_k x^*_3(k) + \tau^h_k x_3) \left(\tau_k/q_k \right)^{0.5-h}$.
According to \eqref{tqxekjx}, at the maximum points $p^*_k$, when $k \to \infty$,
\be
\lab{jxtqxe}
\al
& \f{\tau_k}{q_k} \to S(0), \qq
X_k = \f{\tau_k r^2}{2 q_k}  \to \f{r^2_0}{ 2} S(0), \qq \eta_k   \to x^\infty_3 S^{0.5-h}(0).
\eal
\ee

On the other hand, at this interior maximum point $p^*_k$ of $J_k$, we have
\be
\lab{dx3jk}
\al
&0=\f{\p}{\p x_3} (r^{-1} J_k) =
(0.5+h) \left(\f{\tau_k}{q_k} \right)^{-0.5+h} \underbrace{ \f{\p}{\p x_3} \left(\f{\tau_k}{q_k} \right)}_{O(\tau^{h}_k)} \,  (E + E_2) \\
&\qq + \left(\f{\tau_k}{q_k} \right)^{0.5+h} \f{\p (E+E_2)}{\p X}
 \underbrace{\f{\p}{\p x_3} X_k}_{O(\tau^{h}_k)} +
\left(\f{\tau_k}{q_k} \right)^{0.5+h} \f{\p (E+E_2)}{\p \eta} \underbrace{\f{\p}{\p x_3} \eta_k}_{O(\tau^h_k)}.
\eal
\ee Multiplying both sides by $\tau^{-h}_k$ and using \eqref{p3t/q22}, \eqref{p3Xk2}  and
\eqref{p3etak2} and taking the limit, we find, at the point $p_\infty=(r_0, 0, 0, 0)$, that
\be
\lab{detaEe2}
0=
(0.5+h)  m_1(0) \,  E  + S(0) \f{\p E}{\p X}
 m_2(0) +
S(0) \f{\p E}{\p \eta} m_3(0).
\ee In the above,  $(\f{\p (E+E_2)}{\p X}, \f{\p (E+E_2)}{\p \eta})$ converge on compact sets of $\{X>0 \}$ as $k \to \infty$ by Conditions (a) and (b) of the theorem. Since $p^*_k$ is a maximum point of $J_k$ which converges to the function $J=J(r,  t)$ in $C^0_{x, t, loc}$ sense, the point $p_\infty$  is a maximum point of $J$ whose value is at least $2 \e$.

Substituting \eqref{detaEe2} into \eqref{eqjgudai}, we find $I_2=0$. Hence the left hand side of \eqref{eqjgudai} at $p_\infty$ is less than or equal to $0$ as $p_\infty$ is a maximum point of $J$, while the right hand side is positive i.e.
\[
0 \ge A e^{-A r_0} \left( A +\f{1}{r_0} + \widetilde{R} \right).
\] This contradiction  shows that, for $k$ sufficiently large, \eqref{jk<2e0} holds, i.e.
\be
\lab{jk<2e}
J_k =H_k + \e e^{-A r} \le 2 \e.
\ee

Similarly, by considering the function
\[
\widehat{J}_k \equiv H_k - \e e^{-A r}
\] and working on  its minimal point, one can prove that, for large $k$,
\[
H_k - \e e^{-A r} \ge - 2 \e
\]Therefore
\be
\lab{|hk|<2e}
\Vert H_k \Vert_\infty \le 2 \e.
\ee From \eqref{defhk}, we know that for the scaled variables $(r, x_3, t)$
\be
\lab{defhk2}
\al
H_k(r, \th, x_3, 0)&=r (\tau_k/q_k)^{0.5+h} \left( E(X_k, \eta_k) +  E_2(X_k, \eta_k, \theta, t_k)\right)\\
&=|\tau_k|^{h} \sqrt{x^2_1+x^2_2} \, v_\th(x_1, x_2, x_3, t_k),
\eal
\ee where the last line is under the original physical variable. Since $\e>0$ is arbitrary, by \eqref{|hk|<2e} and \eqref{defhk2}, we deduce
\[
\limsup_{k \to \infty} |\tau_k|^{h} \sqrt{x^2_1+x^2_2} \, |v_\th(x_1, x_2, x_3, t_k)|=0.
\]
This is the final contradiction with \eqref{thrvth2} that proves Theorem \ref{th1.1}. \qed

\section{ Proof of  theorem \ref{th1.2} }

\proof

Suppose $\sig>0.5$. We will reach a contradiction by the end of the section, implying that $\sig=0.5$, i.e. the singularity of $v_r$ must be critical. Then the theorem immediately follows from Theorem \ref{th1.1}.

Recall the divergence free condition from \eqref{eqasns}
\be
\lab{div=0}
 \f{\p (rv_r)}{\p r}+ \f{\p v_\th}{\p \th} +\f{\p ( r v_3)}{\p x_3}=0.
\ee Using the profiles \eqref{sscene2} of the velocity $v$, it turns to
\be
\lab{div=02}
 \f{\p }{\p r} \left( \f{r}{q^\sig} R \right)+ \f{\p }{\p \th} \left( \f{1}{q^{0.5+h}} E_2 \right) +\f{\p }{\p x_3}  \left( \f{r}{q^{0.5+h}} U \right)=0.
\ee Here $R, E_2, U$ are functions of $(X, \eta, \th, t)$. Note the function $E$, the principal angular profile, vanishes since it is independent of $\th$.

As before, let $t_k$ be a sequence such that $t_k \uparrow 0$ and $\tau_k= - t_k$.
Fix any number $x^\infty_3$ and let $x_3(k)= \tau^{0.5-h}_k x^\infty_3$.
Using the variables in \eqref{changvar}, with $x^*_3(k)$ there replaced by $x_3(k)$ here, we can write \eqref{div=02} as
\be
\lab{div=03}
 \f{\p }{\p r} \left( \big(\f{\tau_k}{q_k}\big)^\sig r R \right)+
  \underbrace{\f{\p }{\p \th} \left( \f{\tau^\sig_k}{q^{0.5+h}_k} E_2 \right)}_{I_{3k}} +  \underbrace{\f{\p }{\p x_3}  \left( \f{\tau^\sig_k}{q^{0.5+h}_k} r U \right)}_{I_{4k}}=0.
\ee In the above, $R, E_2, U$ are functions of $(X_k, \eta_k, \th, s_k)$.

According to Condition (d) of the theorem, $\left|\f{\p E_2}{\p \th} \right| \le C r$ where $r$ is the original radial variable. Under the scaled radial variable in \eqref{changvar} line 3, the bound becomes
\[
\left|\f{\p E_2}{\p \th} \right| \le C \tau^{0.5}_k r.
\]Therefore
\be
\lab{i3to0}
|I_{3k}| \le C \tau^{\sig-h}_k (\tau_k/q_k)^{0.5+h} \to 0, \q  k \to \infty.
\ee

Next, notice that
\be
\lab{i4=}
\al
r^{-1} I_{4k} &=  \tau_k^{\sig-0.5-h} \f{\p }{\p x_3}  \left( (\tau_k/q_k )^{0.5+h}  U \right)\\
&= \tau_k^{\sig-0.5-h} \left[ \f{\p }{\p x_3}  \left( (\tau_k/q_k )^{0.5+h} \right)   U
+ (\tau_k/q_k )^{0.5+h} \left( \f{\p U}{\p X} \f{\p X_k}{\p x_3}
+  \f{\p U}{\p \eta} \f{\p \eta_k}{\p x_3} \right) \right] \\
&= \tau_k^{\sig-0.5-h} O(\tau_k^h), \qq \text{by} \q \eqref{p3t/q2}, \eqref{p3Xk}, \eqref{p3etak2}.
\eal
\ee In the above we also used the boundedness of $U, \sqrt{X} \p_X U, \p_\eta U$ (c.f. Conditions (a) and (d) of the theorem). For example, $
| \f{\p U}{\p \eta} \f{\p \eta_k}{\p x_3}| \le C \tau^{h}_k$ and

\[
\al
&\left| \f{\p U}{\p X} \f{\p X_k}{\p x_3} \right| \le C X^{-1/2}_k r^2 \tau^{h}_k
\le C (q_k/\tau_k)^{1/2} r \tau^{h}_k =O(\tau^h_k), \qq \text{by \eqref{changvar} last line}.
\eal
\]
By the contrapositive assumption $\sig>0.5$, we find, on compact sets,
\be
\lab{i4to0}
I_{4k} \to 0, \q k \to \infty.
\ee

Let  $\phi=\phi(r, \th, x_3)$ be a smooth, compactly supported test function. Testing $\phi$ against \eqref{div=03} yields,
\[
\al
0=\int^{2\pi}_0 \int \int \phi (r, \th, x_3) \cdot \left( \f{\p }{\p r} \left( \big(\f{\tau_k}{q_k}\big)^\sig r R \right)+ I_{3k} +  I_{4k} \right)  dr dx_3 d\th.
\eal
\]

Using \eqref{tqxekjx}, \eqref{i3to0} and \eqref{i4to0},  we can take the limit to deduce,
\[
0=\int^{2\pi}_0 \int \int \phi (r, \th, x_3) \cdot  \f{\p }{\p r} \left( r  R(\f{r^2}{2} S(\tau), x^\infty_3 S^{0.5-h}(\tau), \th, 0) \right)
  dr dx_3 d\th.
\]Since $\phi$ is arbitrary, this infers:
for all $r, \tau \ge 0$, $\th \in [0, 2\pi]$,
\be
\lab{RStau}
R(\f{r^2}{2} S(\tau), x^\infty_3 S^{0.5-h}(\tau), \th, 0)=0, \q S(\tau)  = \f{1}{1 +  \tau + S^{-2h}(\tau) (x^\infty_3)^2}.
\ee Taking $\tau=0$, this implies, for $z \equiv  x^\infty_3$,
\be
\lab{eqlz}
L(z) \equiv S(0)+ \left( S^{0.5-h}(0) z \right)^2=1.
\ee We treat the left hand side of this equality as a function of $z \in (-\infty, \infty)$. For any fixed $S(0) \in (0, 1)$, since $L(0)=S(0)<1$ and $\lim_{z \to \infty} L(z) = \infty$, there always exists a solution $z$ for \eqref{eqlz}. Hence  $ S^{0.5-h}(0)  x^\infty_3$ can reach any value in $(0, 1)$ by choosing $x^\infty_3$ suitably. Since $r>0$ is arbitrary, the term $\f{r^2}{2} S(0)$ can reach any positive number. Thus we have proven
\be
\lab{R0=0}
R(X, \eta, \th, 0)=0, \q \forall X>0, \, \eta \in (-1, 1), \, \th \in [0, 2\pi].
\ee
This is a contradiction with Condition (e).
Therefore $\sig=0.5$. As mentioned at the beginning of the proof, this concludes the proof of Theorem \ref{th1.2}. \qed

\begin{remark}
\lab{rmk3.1}
By \eqref{sscene2}, in the original cylindrical system,
\be
\lab{vr=1qsig}
v_r(r, \th, x_3, t) = \frac{1}{q^{\sig}}  R(X, \eta, \th, t).
\ee
The factor $\frac{1}{q^\sig}$  in \eqref{vr=1qsig} can be changed to  functions of the form $\lam= \frac{|\log q|^{\pm 1}}{q^\sig}$, $\sig>0.5$. Indeed, if
$\limsup_{q \to 0} \lam^{-1}(q) |v_r|>0$, then \eqref{R0=0} still leads to a contradiction in exactly the same way.
\end{remark}

\section{Further development and comments}

After the paper was submitted to the arXiv on November 19 2026, two new inputs occurred the next  day, which warrant us to add this short section. 

1. On November 20th, the paper by L. Alpoge and T. Buckmaster " Finite-time blow-up for the Navier-Stokes equations"  appeared on:

https://github.com/tristanbuckmaster/navier-stokes.

The authors stated that under their guidance,   Claude AI was able to prove finite blow up for the unforced Navier-Stokes equations with some smooth compactly supported initial value. Part of the proof is described and lean verification is given.

 Therefore, in the 2nd paragraph in Section 1 of this paper, we should acknowledge this development by adding the above sentences.

\medskip

2. On November 20th, Prof. Xiao Ren kindly asked us to check if the solution in \cite{Op} satisfies the first inequality in Condition (d) in Theorem \ref{th1.2}, i.e.
$\p_\th E_2/r$ is bounded. The answer is no. The reason is that the phase term in the wave function may have derivatives which become unbounded. Therefore Theorem \ref{th1.2} does not cover the solution in \cite{Op}.

\medskip

3. On the other hand, we can modify Conditions (d) and (e)  in Theorem \ref{th1.2} as follows, doing away with the derivative bounds on $\f{1}{r} \p_\th E_2$, $\sqrt{X} \p_X U$, $\p_\eta U$ altogether.

\begin{theorem}
\lab{th4.1}
Let $v=v(x, t)=v_{r}e_{r} + v_{\th}e_{\th} + v_{3}e_{3}$ be a smooth Leray-Hopf solution of the  Navier-Stokes equations in $\mathbf{R^3} \times [-1, 0)$ with a force $f$ that is bounded, smooth and compactly supported in the spatial variables.
Suppose: for $t \in [-1, 0)$,

(a). the functions $E$, $\sqrt{X} \p_X E$, $\p_\eta E$ of variables $(X, \eta)$; $R$  and $U$ of $(X, \eta, \th, t)$  are smooth, bounded  and converge uniformly on compact sets of $\{ X>0 \}$ as $t \uparrow 0$;

(b). the functions
$E_2$, $\sqrt{X} \p_X E_2$, $\p_\eta E_2$  of $(X, \eta, \th, t)$
are smooth, bounded and converge uniformly to $0$ on compact sets of $\{ X>0 \}$ as $t \uparrow 0$; moreover, there are positive numbers $\delta$ and $C$ such that
$|E_2| \le C q^\delta$;

(c).
\[
\sup_{X>0, \eta \in [-1, 1], \th \in [0, 2\pi], t \in [0, 1)} |X|^{0.5^+} (|E(X, \eta)+ E_2(X, \eta, \th, t)|)<\infty
\]where $0.5^+$ is  a number strictly greater than $0.5$;

(d). there is a number $X_0>0$ such that $v_r$ is axially symmetric in $\{X<X_0 \}$;

(e). there is a small number $b>0$ such that $R(\cdot, t)$ is uniformly continuous for $(X, \eta, \th) \in [a, b] \times [-1, 1] \times [0, 2 \pi]$ and $v_r |_{t=0}$ does not vanish on open sets near the $x_3$ axis. Here $b>a>0$.

Then there does not exist a singularity scenario with the following asymptotic profiles
\be
\lab{sscene3}
\al
&v_r = \frac{1}{q^{\sig}}  R(X, \eta, \th, t), \q \sig \ge 0.5,\\
&v_\th = \frac{1}{q^{0.5+h}} \left(E(X, \eta) + E_2(X, \eta, \th, t)\right),\\
&v_3 = \frac{1}{q^{0.5+h}}  U(X, \eta, \th, t), \qq h \in (0, 0.01),
\eal
\ee such that
\be
\lab{thrvth03}
\limsup_{t \to 0^-} |t|^h \Vert r v_\th(\cdot, t) \Vert_\infty>0.
\ee
\end{theorem}

\proof
The proof is similar to that of Theorem \ref{th1.2} except that we pick compactly supported, axially symmetric test functions $\phi=\phi(r, \eta)=\phi(r,  \eta_k(x_3, t))$ in the scaled variables. We choose $\phi$ such that $\p \phi/\p \eta$ is bounded.

Suppose $\sig>0.5$.
Testing $\phi$ against \eqref{div=03} yields, since $I_{3k}=0$ now,
\[
\al
0=\int^{2\pi}_0 \int \int \left( \f{\p }{\p r}  \phi  \cdot  \big(\f{\tau_k}{q_k}\big)^\sig r R  +
\f{\p \phi}{\p \eta} \f{\p \eta_k}{\p x_3}   \f{\tau^\sig_k}{q^{0.5+h}_k} r U \right)  dr dx_3 d\th.
\eal
\] By \eqref{p3etak2},
$
\f{\p}{\p x_3} \eta_k
=O( \tau^{h}_k)$.
Then we can take the limit to deduce, using $\sig>0.5$,
\[
0=\int^{2\pi}_0 \int \int \f{\p }{\p r} \phi (r, x^\infty_3 S^{0.5-h}(\tau)) \cdot   r  R(\f{r^2}{2} S(\tau), x^\infty_3 S^{0.5-h}(\tau), \th, 0)
  dr dx_3 d\th.
\]Alternatively, one can also take $\phi=\phi(r, x_3)$ directly and use \eqref{i4to0} to reach the above equality.

Denoting $\widetilde{R}=  \int^{2 \pi}_0 R d\th/(2\pi)$, this infers
\[
0= \int \int \phi (r, x^\infty_3 S^{0.5-h}(\tau)) \cdot  \f{\p }{\p r} \left( r  \widetilde{R}(\f{r^2}{2} S(\tau), x^\infty_3 S^{0.5-h}(\tau), 0) \right)
  dr dx_3.
\]
As before, we obtain
\be
\lab{R0=02}
\widetilde{R}(X, \eta,  0)=0, \q \forall X>0, \, \eta \in (-1, 1).
\ee
By Condition (d), $R$ is axially symmetric for $X<X_0$, thus
\[
0=\widetilde{R}(X, \eta,  0)= R(X, \eta,  0)
\] for $X <X_b$. Note, by \eqref{taudenj0},  that $\lim_{t  \uparrow 0} q >0$ if $x_3 \neq 0$. By taking $t \uparrow 0$ in line 1 of \eqref{sscene3},  one deduces $v_r |_{t=0}$ vanishes on an open subset of $\{x_3 \neq 0 \}$, contradicting Condition (e).
Therefore $\sig=0.5$. By Theorem \ref{th1.1}, this  concludes the proof of Theorem \ref{th4.1}. \qed

\medskip

The solution in \cite{Op} satisfies Condition (d) in this theorem by construction in Sections 4-5. Note there is only one purported singular point at $t=0$, the origin.
It  also satisfies Conditions  (e) due to the almost explicit construction using analyticity near the axis.

Using the notations in \eqref{sscene3} and \eqref{opai1} (solution in \cite{Op}),
\be
\lab{e2zhankai}
q^{-(0.5+h)} E_2 = u^{ap} \cdot e_\th= w \cdot e_\th + \sum u^{[i]} \cdot e_\th= w \cdot e_\th + \text{lower order terms}.
\ee Here "lower order terms" stand for terms with higher powers of $q$. According to p71 of \cite{Op}, the  wave function $w$ can be written as
\be
\lab{wzhankai}
w= \underbrace{\left( a^{(r)} e_r + a^{(\th)} e_\th + a^{(3)} e_3 \right)}_{AM} \exp [i (N \th + g(r, x_3))]
+  \text{lower order terms}.
\ee In the above $N$ is a finite integer, $g$ is a $C^1$ function and $a^{(r)}, a^{(\th)}$ and $a^{(3)}$ are amplitude functions independent of $\th$. Eventually one also needs to take the real part. Here we have simplified the notations a little. In fact $N$ is   determined by $k_{rv} m_v p_r$ defined in (6.27), with $k_{rv}$ being  $p^{-1}_r$ essentially. The amplitude of $|w|$ is in the order of $q^{-0.5-0.5 h}$
\be
\lab{pthe2}
E_2 = q^{0.5 + h} (w \cdot e_\th)| = O(q^{0.5 h}) \exp [i (N \th + g(r, x_3))] +  \text{lower order terms}.
\ee Due to the term $q^{0.5 h}$,  it seems plausible that Condition (b) holds, although we have not finished checking.

To recap, to check if the OpenAi solution is covered by Theorem \ref{th4.1}, the main remaining task is to examine if  the condition on $\sqrt{X} \p_X E_2$, $\p_\eta E_2$ in Condition (b) hold in an annulus $0<X_a<X<X_b$.

\section*{Acknowledgments} We wish to thank Professors  Igor Kukavica,  Zhen Lei, Zijin Li,
  Xiao Ren, Xin Yang, Qi Xu for helpful discussions.

\centerline {\it Statements and Declarations.}

 1. The paper is written in the traditional way without AI use in any form. It is likely to be the last paper written by the author in this manner.  2. There is no data involved in this paper.

\bibliographystyle{plain}



\def\cprime{$'$}

\end{document}